\documentclass[review, 3p]{elsarticle}
\makeatletter
\def\ps@pprintTitle{%
	\let\@oddhead\@empty
	\let\@evenhead\@empty
	\def\@oddfoot{\centerline{\thepage}}%
	\let\@evenfoot\@oddfoot}
\makeatother
\usepackage{amsmath}
\usepackage{amssymb}
\usepackage{mathtools}
\usepackage{enumerate}

\usepackage{multicol}
\usepackage{todonotes}

\usepackage{hyperref}
\numberwithin{equation}{subsection}

\usepackage{amsthm}
\theoremstyle{plain}
\newtheorem{theorem}{Theorem}[section]

\newtheorem{lemma}[theorem]{Lemma}
\newtheorem{proposition}[theorem]{Proposition}
\newtheorem{corollary}[theorem]{Corollary}
\theoremstyle{remark}
\newtheorem*{remark}{Remark}

\theoremstyle{definition}
\newtheorem{definition}[theorem]{Definition}

\newtheorem*{notation}{Notation}

\makeatletter
\@ifpackageloaded{hyperref}%
{\newcommand{\mylabel}[2]% #1=name, #2 = contents
	{\protected@write\@auxout{}{\string\newlabel{#1}{{#2}{\thepage}%
				{\@currentlabelname}{\@currentHref}{}}}}}%
{\newcommand{\mylabel}[2]% #1=name, #2 = contents
	{\protected@write\@auxout{}{\string\newlabel{#1}{{#2}{\thepage}}}}}
\makeatother

\newcommand{\F}{\mathbb{F}}
\newcommand{\Prob}{\mathbb{P}}
\newcommand{\PI}{\mathbb{P}_{\textrm{inv}}}

\renewcommand{\epsilon}{\varepsilon}

\newcommand{\tr}{\textrm{tr}}
\newcommand{\Inv}{\textrm{Inv}}
\newcommand{\Stab}{\textrm{Stab}}
\newcommand{\rad}{\textrm{rad}}
\newcommand{\diag}{\textrm{diag}}
\newcommand{\B}[1]{\mathbb{#1}}
\newcommand{\Property}{\mathcal{P}}
\newcommand{\I}{\mathcal{I}}
\newcommand{\T}{\mathcal{T}}

\begin{document}

	\begin{frontmatter}
	\title{Invariable generation of finite classical groups}
	\author{Eilidh McKemmie\fnref{myfootnote}}
	\fntext[myfootnote]{This work was taken from the author's PhD thesis. The author thanks her PhD advisor Robert Guralnick and the reviewer for their time and very helpful comments and acknowledges the support from the NSF under grant DMS 1600056.}
	\address{Institute of Mathematics, Hebrew University, Jerusalem, Israel}
	\ead{eilidh.mckemmie@gmail.com}
	
	\begin{abstract}
		A subset of a group invariably generates the group if it generates even when we replace the elements by any of their conjugates. In a 2016 paper, Pemantle, Peres and Rivin show that the probability that four randomly selected elements invariably generate $S_n$ is bounded away from zero by an absolute constant for all $n$. Subsequently, Eberhard, Ford and Green have shown that the probability that three randomly selected elements invariably generate $S_n$ tends to zero as $n \rightarrow \infty$. In this paper, we prove an analogous result for the finite classical groups. More precisely, let $G_r(q)$ be a finite classical group of rank $r$ over $\F_q$. We show that for $q$ large enough, the probability that four randomly selected elements invariably generate $G_r(q)$ is bounded away from zero by an absolute constant for all $r$, and for three elements the probability tends to zero as $q \rightarrow \infty$ and $r \rightarrow \infty$. We use the fact that most elements in $G_r(q)$ are separable and the well-known correspondence between classes of maximal tori containing separable elements in classical groups and conjugacy classes in their Weyl groups.
	\end{abstract}

\end{frontmatter}
	
	\section{Introduction}
	
	Motivated by a problem in computational Galois theory, Dixon \cite{Dixon2} began studying random invariable generation of the symmetric group in 1992.
	
	\begin{definition}
		Let $G$ be a finite group and let $S=\{s_1, ..., s_k\}\subseteq G$. The subgroup of $G$ invariably generated by $S$ is \[\bigcap_{g_1, ..., g_k \in G}\langle s_1^{g_1}, ..., s_k^{g_k}\rangle\]
	\end{definition}

	Note that $S$ invariably generates $G$ if and only if $S$ is a generating set that remains a generating set upon replacing any element with a conjugate.

	Consider the problem of deciding whether or not the Galois group of a polynomial $f(x) \in \B{Z}[x]$ of degree $n$ is $S_n$. As highlighted by Pemantle, Peres and Rivin \cite{Pemantle}, there are deterministic algorithms with complexity $O(n^{40})$, but typically the running time is too high to be practical. There is also a Monte Carlo algorithm which relies on the Frobenius Density Theorem. Recall that we may factorise $f$ modulo a prime $p$ into distinct irreducible factors, except in the finitely many cases where $p \mid \textrm{disc}(f)$. Let $\lambda_p$ be the partition of $n$ consisting of the degrees of these factors.
	
	For a partition $\lambda$ of $n$, the Frobenius Density Theorem (see \cite[IV.5.2]{janusz1996algebraic} for example) states that the density of primes $p \nmid \textrm{disc}(f) $ for which $\lambda_p=\lambda$ is equal to the density of elements in $Gal(f)$ which correspond to the partition $\lambda$. Thanks to this theorem, even though we do not know $Gal(f)$, we may uniformly sample $k$ random elements and get their cycle types by randomly picking primes $p_1, ..., p_k$ and finding the $\lambda_{p_i}$. If elements of these cycle types invariably generate $S_n$ then the algorithm tells us that $Gal(f)=S_n$, which we know to be true. If they do not invariably generate then the algorithm tells us that $Gal(f)\neq S_n$ although this is not always accurate. We get a false negative if we have picked $k$ random elements in $S_n$ which do not invariably generate. To understand how useful this algorithm is, we would like to know the probability that $k$ random elements in $S_n$ invariably generate.
	
	\begin{definition}
		For a finite group $G$, a positive integer $k$, a subset $H \subseteq G$ and a set $\mathcal{A}$ of subgroups of $G$, define
		\[\PI(G,k,H,\mathcal{A})=\frac{|\{(g_1, ..., g_k) \in H^k \mid \{g_1, ..., g_k\} \mbox{ invariably generates a subgroup in }\mathcal{A}\}|}{|G^k|}\] which is the proportion of elements of $G^k$ which are $k$-tuples of elements of $H$ and invariably generate a member of $\mathcal{A}$. For simplicity, we will omit $H$ when $H=G$ and we will omit $\mathcal{A}$ when $\mathcal{A}=\{G\}$. So $\PI(G,k)$ is the probability that $k$ random elements of $G$ invariably generate $G$.
	\end{definition}

	If $G$ is a permutation group, let $\mathcal{T}$ be the set of proper transitive subgroups of $G$. The following theorem on random invariable generation of $A_n$ and $S_n$ is proved in \cite{Pemantle,EFG2}.
	
	\begin{theorem}\mylabel{efg}{Theorem~\thetheorem{}}
		Let $G$ be $A_n$ or $S_n$.		
		\begin{enumerate}[(a)]
			\item There is a constant $b>0$ such that $\PI(G,4,\mathcal{T})\ge \PI(G,4)\ge b$ for all $n$.
			\item For three elements, $\lim_{n \rightarrow \infty}\PI(G,3,\mathcal{T})=\lim_{n \rightarrow \infty}\PI(G,3)=0$.
		\end{enumerate}
	\end{theorem}

	This allows us to determine the accuracy of the Monte Carlo algorithm described above. The probability of getting a false negative after picking $k$ random primes is at most $(1-b)^{k \choose 4}$.

	The approach in \cite{Pemantle,EFG2} involves modelling the number of $j$-cycles of a random permutation as a random Poisson variable with mean $\frac{1}{j}$ (which is discussed by Arratia and Tavar{\'e} \cite{Arratia_1992} for example) and noting that a subset of $S_n$ fails to invariably generate a transitive subgroup of $S_n$ if and only if the elements all fix a set of size $j$ for some $1\le j <n$. Eberhard, Ford and Green show that a typical permutation does not have many small cycles, and that four permutations with few small cycles are unlikely to fix a set of a common size. Using a result of {\L}uczak and Pyber \cite{pyber} which states that the union of proper transitive subgroups of $S_n$ not containing $A_n$ is small, they are able to obtain the desired result for $A_n$ and $S_n$. For three elements, they take inspiration from Maier and Tenenbaum's \cite{maier1984set} study of the distribution of prime factors of random integers. They show that the number of common fixed set sizes of three random permutations grows unboundedly with $n$ with high probability.
	
	More generally, there has been recent interest in invariable generation properties of finite simple groups. It is well known \cite{Aschbacher_1984, Steinberg} that every finite simple group is generated by two elements. In 1882, Netto \cite{Netto} conjectured that a random pair of elements in the symmetric group generates with probability tending to $\frac{3}{4}$ as the size of the group approaches infinity. Dixon \cite{Dixon} proved the conjecture in 1969, reigniting activity in probabilistic group theory. Dixon's work was continued by Kantor and Lubotzky \cite{KL} and Liebeck and Shalev \cite{LS} to show that two random elements generate a finite simple group $G$ with probability approaching $1$ as $|G| \rightarrow \infty$.
	
	Are there invariable versions of these theorems? Guralnick and Malle \cite{GM}, and independently Kantor, Lubotzky and Shalev \cite{Kantor2011} show that every finite simple group is invariably generated by two elements, but it is not true that almost all pairs invariably generate. Indeed, Kantor, Lubotzky and Shalev \cite{Kantor2011} showed that there is an absolute constant $\epsilon>0$ such that $\PI(G,k)\le (1-\epsilon)^k$ for all finite simple groups $G$ and positive integers $k$. There are also some results for infinite groups. For instance, Kantor, Lubotzky and Shalev \cite{Kantor2015} showed that a linear group is invariably generated by a finite set if and only if it is finitely generated and virtually solvable.
	
	The aim of this paper is to begin extending the results in \ref{efg} to the finite simple groups. In particular, we want to know how many random elements we need to pick in order to bound the probability of invariable generation away from zero by a constant that does not depend on the order of the group.
	
	\begin{notation}
		Let $K=\overline{\F_p}$ be the algebraic closure of the field $\F_p$ of positive characteristic $p$, and let $X$ be a simply connected linear algebraic group (not necessarily simple) of rank $r$ over $K$ acting on a $K$-vector space $E$. We will consider $X=SL_n(K), Sp_{2n}(K)$ with the natural matrix representation and $X=Spin_{m}(K)$ (assuming that $p$ is odd when $m$ is odd) acting on $E:=K^m$ via the surjection \begin{equation}\label{surjection}
		\phi:Spin_m(K)\rightarrow SO_m(K)
		\end{equation}
		which has a kernel of order $2$ in odd characteristic and is an isomorphism in even characteristic.
		
		Denote by $W$ the Weyl group of $X$. Let $\sigma:X \rightarrow X$ be a Steinberg endomorphism such that the set of fixed points $X_\sigma$ is a finite group of Lie type over $\F_q$ where $q$ is a power of $p$. We refer to the groups in the following list as the finite classical groups: \begin{equation}\label{classical groups}SL_n(q) \mbox{ and } Sp_{2n}(q) \mbox{ for }n\ge2, SU_n(q) \mbox{ and } \Omega^\pm_{2n}(q) \mbox{ for }n \ge 3, \Omega_{2n+1}(q) \mbox{ for }q\mbox{ odd and }n\ge 2\end{equation} and note that for $n\ge 4$ these groups are quasisimple (in fact, for $n=3$ all of these groups except $SU_3(2)$ are quasisimple). Every finite classical group $G$ is of the form $X_\sigma$ or a quotient of $X_\sigma$ by a central subgroup. We will abuse notation by viewing $G$ as a subgroup of $X$ when really it is the lift of $G$ to $X_\sigma$ that is a subgroup of $X$, but this should not cause any confusion.
	\end{notation}
	
	In joint work with Daniele Garzoni \cite{Garzoni}, we fix $r$ and let $q \rightarrow \infty$, and show that two random elements of $G_r(q)$ invariably generate with probability bounded away from zero.
	
	\begin{theorem}\mylabel{GM}{Theorem~\thetheorem{}}
		There are constants $\epsilon_r>0$ depending only on $r$ and not on $q$ such that if $G_r(q)$ is a finite group of Lie type then $\PI(G_r(q),2)\ge \epsilon_r$.
	\end{theorem}

	We give an explicit formula for estimating $\PI(G_r(q),2)$ with an error term of $O(1/q)$ for fixed $r$ \cite[Theorem~6.2]{Garzoni}. The values are quite easy to compute, for example $\PI(G_2(3^a),2)=1/9+O(1/q)$ and $\PI(SL_2(q), 2)=1/2+O(1/q).$

	In this paper we will let $r \rightarrow \infty$. Our main result is the following.
	
	\begin{theorem}\mylabel{main}{Theorem~\thetheorem{}}
		Let $G=G_r(q)$ be a finite group of Lie type of rank $r$.
		\begin{enumerate}[(a)]
			\item There exist absolute constants $c\in (0,1)$ and $Q \in \mathbb{N}$ such that for all $q > Q$ and for all $r$ we have $\PI(G,4)>c$.
			\item $\lim_{q \rightarrow \infty}\lim_{r \rightarrow \infty}\PI(G,3)=0$.
		\end{enumerate}
	\end{theorem}
	
	The exceptional groups have bounded rank and so for part $(a)$ they are covered by \ref{GM} and we do not need to consider them when proving $(b)$. Therefore we will work with $G$ a finite classical group from the list (\ref{classical groups}) until the final section of this paper. Our proof relies in a crucial way on \ref{efg} and the proof of Eberhard, Ford and Green. In the first five sections of this paper we will show that $\PI(G,k)$ is close to $\PI(S_n,k)$ for some integer $n$. In Section~\ref{separable} we will define separable elements to be those with distinct eigenvalues, and show that it is enough to study the probability that random separable elements invariably generate an irreducible subgroup of $G$. Section~\ref{invariant} describes an action of $G$ on a module $V$ and an action of the Weyl group on a set $M$, and observes a relationship between invariant subspaces of $V$ and invariant subsets of $M$. In Section~\ref{quasi} we show it is enough to know the dimensions and types of invariant subspaces of random separable elements of $G$, and in Section~\ref{irred sep} we show that this is roughly the same as understanding the fixed sets with respect to a natural action of a random element of the Weyl group $W$. In Section~\ref{weyl} we give a modified version of Eberhard, Ford and Green's proof of \ref{efg}(b) which allows our argument to work for the orthogonal groups. Finally in Section~\ref{proof} we complete the proof of \ref{main} for finite classical groups and extend the result to the other finite groups of Lie type.
	
	\section{Separable elements and irreducible subgroups}\label{separable}
	
	Observe that a subset $S$ of a finite classical group $G$ fails to invariably generate an irreducible subgroup if and only if each element can be conjugated to stabilise a common subspace of the natural module. Therefore we will be interested in invariant spaces of elements in $G$. We will focus on the set of separable elements since their invariant spaces are easy to describe and a large proportion of elements are separable.
	
	Let $G$ be a finite classical group from the list (\ref{classical groups}). The ambient simply-connected algebraic group $X$ acts on the $K$-vector space $E$ -- in the case of the spin groups it acts via the surjection $\phi$ defined in (\ref{surjection}) -- so that elements of $X$ may be written as matrices with entries in $K$. An element $s \in X$ is called separable if it has no multiple eigenvalues. Denote by $s(G)$ the set of separable elements of $G$ (that is, the elements in $G$ which lift to separable elements of $X$).
	
	If $s\in X$ is separable then $C_X(s)$ is the unique maximal torus of $X$ containing $s$ (see, for example, \cite[Corollary~14.10]{MalleTesterman}). When $s \in X_\sigma$, this torus is $\sigma$-stable. In every case, there exist separable elements in $G$. In fact, it was shown in {\cite[Theorem~1.1]{GuLub}} and {\cite[Theorems~6.1,~6.2]{FulmanGuralnick}} that most elements of $G$ are separable. Define \begin{equation}\label{epsilon}\epsilon(q)=\frac{7}{q-1}-\frac{2}{(q-1)^2}.\end{equation}
	
	\begin{proposition}\mylabel{most sep}{Proposition~\thetheorem{}}
		Let $G$ be a finite classical group from the list (\ref{classical groups}) over $\F_q$. Then
		\[
		\frac{\left |s(G)\right |}{\left| G \right |}\ge 1-\epsilon(q).
		\]
	\end{proposition}
	
	Let $G$ act on its natural module $V$ and let $\I$ be the set of irreducible subgroups of $G$ (note that $G\in \I$). If $q$ is even then, as we will see in Section~\ref{symplectic section}, every separable element of $Sp_{2n}(q)$ is contained in a subgroup isomorphic to $O_{2n}^+(q)$ or $O_{2n}^-(q)$ but not both. Define
	\[\I'=\begin{cases}\{H \in \I \mid H\ncong O^\pm_{2n}(q)\} &
	G=Sp_{2n}(q), q \mbox{ even}\\ \I & \mbox{otherwise}.\end{cases}\]
	
	Let $Y(G)$ be the set of separable elements contained in any element of $\I'\setminus \{G\}$.
	
	\begin{proposition}\mylabel{rss irred}{Proposition~\thetheorem{}}
		There is a function $\alpha:\mathbb{N}\times \mathbb{N} \rightarrow [0,1]$ such that $\lim_{n\rightarrow \infty}\alpha(n, k)=0$ for any $k$ and
	\begin{equation*}\PI(G, k, s(G), \I') - \alpha(n, k) \le \PI(G,k)\le \PI(G, k, s(G), \mathcal{I})+\epsilon(q)^k.\end{equation*}
\end{proposition}
\begin{proof}
	By definition $\PI(G, k, s(G)) \le \PI(G, k)$ and
	\[
	\PI(G, k)\le \PI(G, k, s(G)) + \left(1-\frac{\left |s(G)\right |}{\left| G \right |}\right)^k.
	\]
	By \ref{most sep} \[\PI(G, k)\le \PI(G, k, s(G)) + \epsilon(q)^k\] and by definition, $\PI(G, k, s(G))\le \PI(G, k, s(G), \I)$ and
	\[
	\PI(G, k, s(G), \I')-\alpha(n, k) \le \PI(G, k, s(G)) \quad \mbox{ where }\alpha(n, k):=\left(\frac{|Y(G)|}{|G|}\right)^k.
	\] Therefore $\PI(G, k, s(G), \I') - \alpha(n, k) \le \PI(G,k)$.
	
	Finally, \cite[Theorem~1.7]{FulmanGuralnick4} implies that $\lim_{n\rightarrow \infty}\alpha(n, k)=0$ for all $k$ and $q$.
\end{proof}
	
	Now we will concentrate on the probability that randomly chosen separable elements invariably generate an irreducible subgroup of $G$.
	
	The set $T=R_\sigma$ of $\sigma$-fixed points of a $\sigma$-stable maximal torus $R$ of $X$ will be called a maximal torus of $G$. Every separable element in $G$ is contained in a unique maximal torus.
	
	\begin{proposition}\mylabel{torus element}{Proposition~\thetheorem{}}
		Let $G$ be a finite classical group from the list (\ref{classical groups}). If $g \in G$ is separable and contained in a maximal torus $T$ of $G$, then $g$ and $T$ have the same invariant subspaces on the natural module $V$.
	\end{proposition}
	\begin{proof}
		Clearly, any space stabilised by $T$ is also stabilised by $g$. Let $\Stab_X(U)$ be the stabiliser in $X$ of a space $U<V$. If $g\in \Stab_X(U)$ then $U$ is spanned by eigenvectors of $g$ and since $g$ is separable, the maximal torus $C_X(g)$ of $X$ stabilises $U$. Since separable elements are contained in a unique maximal torus, it must be that $T =C_X(g)\le \Stab_X(U)$.
	\end{proof}
	
	The Steinberg endomorphism $\sigma$ induces an automorphism of the Weyl group which we will also refer to as $\sigma$, and we consider the semidirect product $W.\langle \sigma \rangle$, and the coset $W\sigma$. If $\sigma$ acts trivially on $W$ then $W.\langle \sigma \rangle \cong W$.
	
	We will exploit the following correspondence in order to apply \ref{efg} in the proof of our main result. For a proof of the correspondence, see \cite[Proposition~25.1]{MalleTesterman}, for example.
		
	\begin{proposition}
		\mylabel{correspondence}{Proposition~\thetheorem{}}
		There is a bijection
		\[\{G\mbox{-classes of }\sigma\mbox{-stable maximal tori of }X\} \leftrightarrow \{W\mbox{-conjugacy classes in }W\sigma\}.\]
	\end{proposition}

	Specifically, fix a $\sigma$-stable maximal torus $R$ of $X$ so that $N_X(R)/R \cong W$, and let $x \in X$ be such that $w = x^{-1}\sigma(x) R\in N_X(R)/R \cong W$. Then $R_w:=R^x$ is the $\sigma$-stable maximal torus corresponding to the $W$-class of $w.\sigma$. If $\rho:X_\sigma \rightarrow G$ is the quotient map then we write $T:=\rho\left(R_\sigma\right)$ and $T_w := \rho\left(\left(R_w\right)_{\sigma}\right)$ where \begin{equation}\label{compute Tw}\left(R_w\right)_{\sigma} \cong \left\{t \in R \mid t = \sigma\left(t\right)^{x^{-1}\sigma(x)}\right\}.\end{equation}
	
	Every separable element $s \in G$ is contained in a unique maximal torus $C_X(s)$ of $X$  which is conjugate to some $R_w$. This allows us to associate to $s$ the $W$-conjugacy class $(w.\sigma)^W$ of $W.\langle \sigma \rangle$.
	
	\section{Tori and their invariant spaces}
	\label{invariant}
	
	In this section, we define a set $M$ acted on by the Weyl group $W$, determine the invariant sets of $w\in W$ on $M$, and use \ref{correspondence} to describe the invariant subspaces of $T_w\le X_\sigma$ in its action on $V$.
	
	When $V$ is a vector space with a bilinear form $f$, define the radical of a subspace $U$ to be $\rad(U)=U \cap U^{\perp}$. Call $f$ non-singular if $\rad(f):=\rad(V)=0$. A space $U$ is called non-singular if $\rad(U)=0$ and totally singular if $\rad(U)=U$. When $V$ is a vector space with a quadratic form $Q$, there is an induced bilinear form $f(u,v)=Q(u+v)-Q(u)-Q(v)$. Define $\rad(Q)=\{v \in \rad(f) \mid Q(v)=0\}$ and call $Q$ non-degenerate if $\rad(f)=0$ and non-singular if $\rad(Q)=0$. Unless $V$ is odd-dimensional over a field of even characteristic, $Q$ is non-degenerate if and only if it is non-singular. Each even-dimensional non-singular subspace of $V$ has either plus or minus type: if $U<V$ has dimension $2m$ then it is of plus type if there is a totally singular subspace of dimension $m$, and minus type otherwise. When $V$ has odd dimension and a non-singular quadratic form, we also refer to the type of a proper non-singular odd-dimensional subspace $U$ by which we mean the type of the orthogonal complement $U^\perp$.
	
	\begin{definition}\mylabel{spacesdef}{Definition~\thetheorem{}}
	Let $G$ be a finite classical group from the list in (\ref{classical groups}). For $g \in G$ and $\Property$ some set of non-zero proper subspaces of $V$ define \[\Inv_\Property(g):=\{\dim(U) \mid U \in \Property, g(U)=U\},\] the set of dimensions of $g$-invariant spaces in $\Property$. If $\Property$ is the set of all non-zero proper subspaces of $V$ we omit the subscript $\Property$.
	\end{definition}

	We will use the same language for analogous notions in the Weyl group. Denote by $B_n$ the hyperoctahedral group $C_2 \wr S_n\le S_{2n}$. The Weyl group $W$ of $X$ is contained in $B_n$ for some $n$, so we will describe some useful properties of this group. We will often think of $B_n$ as the group of permutations $w$ of $\{\pm 1, ..., \pm n\}$ with the property that $w(-i)=-w(i)$. By \cite[Section~4.2]{JK}, the conjugacy classes of $B_n$ correspond to signed cycle types, that is, $(n_1^{\epsilon_1}, \dots, n_r^{\epsilon_r})$ where $\epsilon_i=\pm1$ and $n_1, \dots, n_r$ are positive integers which sum to $n$. A positive $j$-cycle is $(a_1, ..., a_j)(-a_1, ..., -a_j)$ for some $a_i \in \{\pm 1, ..., \pm n\}$ such that the $|a_i|$ are distinct, and a negative $j$-cycle is $(a_1, ..., a_j, -a_1, ..., -a_j)$ for some $a_i \in \{\pm 1, ..., \pm n\}$ such that the $|a_i|$ are distinct. The sign of an element with signed cycle type $(n_1^{\epsilon_1}, \dots, n_r^{\epsilon_r})$ is $\prod_i \epsilon_i$ and there is a subgroup $B_n^+\le B_n$ of index $2$ consisting of the positive elements of $B_n$.
	
	\begin{definition}
		\begin{enumerate}
			\item For $w \in W$ acting on a set $M$ and $\Property$ some set of non-empty proper subsets of $M$ define \[\Inv_\Property(w):=\{|S| : S\in \Property, w(S)=S\},\] the set of sizes of the $w$-invariant sets in $\Property$. If $\Property$ is the set of all non-empty proper subsets of $M$ we omit the subscript $\Property$.
			\item When $M\subseteq\{\pm 1, ..., \pm n\}$, we call a subset $S \subseteq M$ totally singular if whenever $k \in S$ then $-k \notin S$, and non-singular when $k \in S$ if and only if $-k \in S$.% In particular every subset of $M=\{1, ..., n\}$ is totally singular.
			\item An element in the conjugacy class corresponding to $(n_1^{\epsilon_1}, \dots, n_r^{\epsilon_r})$ is said to be of plus type if $\prod \epsilon_i=1$ and of minus type if $\prod \epsilon_i=-1$. If $S\subseteq \{\pm 1, ..., \pm n\}$ is non-singular and stabilised by $w$ then we may restrict $w$ to $S$ and define the type of $S$ under $w$ to be the type of $w|_S$.
		\end{enumerate}
	\end{definition}

	For example, consider the action of $B_4$ on $M=\{\pm1, ..., \pm4\}$ and take $w=(1, 2, 3)(-1, -2, -3)(4, -4)$ which has signed cycle type $(3^+, 1^-)$. The non-singular sets stabilised by $w$ are \[S_1 = \{1, 2, 3, -1, -2, -3\}, \quad
	S_2 = \{4, -4\}, \quad S_1 \cup S_2\] and under $w$ the set $S_1$ is of plus type while $S_2$ and $S_1 \cup S_2$ are of minus type.
	
	\begin{definition}
		\begin{enumerate}
			\item By $\mathcal{N}$ denote both the set of non-zero proper non-singular subspaces of $V$ and the set of non-empty proper non-singular subsets of $M$.
			\item Denote by $\mathcal{R}$ both the set of non-zero proper totally singular subspaces of $V$ and the set of non-empty proper totally singular subsets of $M$.
			\item Define $\mathcal{N}^+=\{N \in \mathcal{N} \mid N\mbox{ has plus type}\}$ and $\mathcal{N}^-=\{N \in \mathcal{N} \mid N\mbox{ has minus type}\}$.
		\end{enumerate}
	This should cause no confusion because it is clear when we are working in $V$ or $M$.
	\end{definition}
		
	\subsection{Linear groups}
		Let $X=SL_n(K)$ with natural module $E$ and take the Steinberg endomorphism $\sigma:(a_{ij}) \mapsto (a_{ij}^q)$ so that $G=SL_n(q)$. Let $V$ be the natural module of $G$ with basis $\{v_1, ..., v_n\}$ and, for each $d$, note that $G$ acts transitively on the set of subspaces of $V$ of dimension $d$. Taking our maximal torus $R$ to be the group of diagonal matrices in $X$, the Weyl group is the set of permutation matrices $W \cong S_n$ upon which $\sigma$ acts as the identity. Therefore $W.\langle \sigma \rangle\cong W$. We work with the natural action of $W$ on the set $M=\{1, ..., n\}$. By \ref{correspondence} the $G$-classes of $\sigma$-stable maximal tori of $X$ correspond to conjugacy classes of $S_n$ which correspond to partitions of $n$.
		
		\begin{proposition}\mylabel{linear}{Proposition~\thetheorem{}}
			If $w\in W$ is associated to the partition $n=n_1+ \cdots + n_m$ then $V$ decomposes as $V=V_1 \oplus \cdots \oplus V_m$ where $T_w$ acts irreducibly on $V_i$ and $\dim V_i = n_i$ for all $i$.
		\end{proposition}
		\begin{proof}
			By equation (\ref{compute Tw}) in \ref{correspondence}, $T_w \cong D_1 \times \cdots \times D_m$ where, written over the algebraic closure $K$, \[D_i\cong \left\{\diag\left(t, t^q, ..., t^{q^{n_i-1}}\right) : t^{q^{n_i}-1}=1\right\}\] which acts irreducibly on the space $V_i$ with basis $\{v_j \mid \sum_{k=1}^{i-1}n_k < j \le \sum_{k=1}^{i}n_k\}$.
		\end{proof}
	\subsection{Symplectic groups}\label{symplectic section}
	
	Let $X=Sp_{2n}(K)$ acting on the $K$-vector space $E$ with basis $\{v_1, ..., v_{n}, v_{-n}, ..., v_{-1}\}$ preserving a non-singular alternating bilinear form \[\left(\, , \,\right):E \times E \rightarrow K, \quad \left( \sum_{i=-n}^n x_iv_i,\sum_{i=-n}^n y_iv_i\right) = \sum_{i=-n}^n x_iy_{-i}.\] Take the Steinberg endomorphism $\sigma:(a_{ij}) \mapsto (a_{ij}^q)$. Then $G=Sp_{2n}(q)$ which acts on the vector space $V$ over $\F_q$ with basis $\{v_1, ..., v_{n}, v_{-n}, ..., v_{-1}\}$, and preserves the restriction of the form $\left(\, , \,\right)$ to $V$. For each $d$, there is exactly one $G$-orbit of non-singular spaces of dimension $d$, and exactly one $G$-orbit of totally singular spaces of dimension $d$. Taking the maximal torus $R:=\{\diag(a_1, ..., a_n, a_{n}^{-1}, ..., a_{1}^{-1}) \mid a_i \in K^\times\}$ of $X$, we can write the Weyl group as a set of permutation matrices in $S_{2n}$ isomorphic to $B_n$. Therefore $\sigma$ acts trivially on $W$ and so $W.\langle\sigma\rangle\cong W\cong B_n$. We will let $W$ act on $M:=\{\pm 1, ..., \pm n\}$ in the natural way. The $G$-classes of $\sigma$-stable maximal tori of $X$ correspond to conjugacy classes of $B_n$ which correspond to signed partitions of $n$.
	
	\begin{proposition}\mylabel{symplectic}{Proposition~\thetheorem{}}
		If $w \in W$ is in the conjugacy class associated to $(n_1^{\epsilon_1}, \dots, n_m^{\epsilon_m})$ then $V$ decomposes under the action of $T_w$ as $V=V_1 \oplus \cdots \oplus V_m$ such that $\dim V_i=2n_i$ for each $i$. If $\epsilon_i=1$ then $V_i$ is a non-singular sum of two irreducible totally singular $n_i$-spaces, whereas if $\epsilon_i=-1$ then $V_i$ is a non-singular space upon which $T_w$ acts irreducibly.
	\end{proposition}
	\begin{proof}
		In order to find $T_w$, we take an element $r=\diag(a_1, ..., a_n, a_{n}^{-1}, ..., a_{1}^{-1}) \in R$ and let $x^{-1}\sigma(x)$ be the permutation matrix representing $w$. By (\ref{compute Tw}) we see that if $r \in T_w$ then $r=\sigma(r)^{x^{-1}\sigma(x)}$. If $\epsilon_1=1$ then the computation gives us $(a_1, ..., a_{n_1})=(a_2^q, ..., a_{n_1}^q, a_1^q)$, whereas if $\epsilon_1=-1$ then we get $(a_1, ..., a_{n_1})=(a_2^q, ..., a_{n_1}^q, a_1^{-q})$, and each cycle gives us similar results. So $T_w\cong D_1^{\epsilon_1} \times \cdots \times D_m^{\epsilon_m}$ where, over the field $K$, the group $D_i^+$ may be written\[D_i^+\cong \left\{\diag\left(t, t^q, ..., t^{q^{n_i-1}}, t^{-q^{n_i-1}},..., t^{-q}, t^{-1}\right) : t^{q^{n_i}-1}=1\right\}\] and $D_i^-$ may be written
		
		\[D_i^-\cong \left\{\diag\left(t, t^q, ..., t^{q^{n_i-1}}, t^{-q^{n_i-1}},..., t^{-q}, t^{-1}\right) : t^{q^{n_i}+1}=1\right\}.\]
		
		We see that $D_i^{\epsilon_i}$ stabilises the space $V_i=V_i^+ \oplus V_i^-$ where $V_i^+$ has basis $\{v_{j} \mid \sum_{k=1}^{i-1}n_k < j \le \sum_{k=1}^{i}n_k\}$ and $V_i^-$ has basis $\{v_{-j} \mid \sum_{k=1}^{i-1}n_k < j \le \sum_{k=1}^{i}n_k\}$.  The group $D_i^+$ acts irreducibly on the totally singular spaces $V_i^+, V_i^-$ while $D_i^-$ acts irreducibly on the non-singular space $V_i$.
	\end{proof}

	In characteristic $2$ we have $Sp_{2n}(q)\cong O_{2n+1}(q)$ which preserves a symmetric bilinear form $f$ with a $1$-dimensional radical. The group acts trivially on the radical so that each element preserves a bilinear symmetric form on the $2n$-dimensional space $\rad(f)^\perp$. So every element of $Sp_{2n}(q)$ is contained in a subgroup isomorphic to $O^+_{2n}(q)$ or $O^-_{2n}(q)$. By \ref{symplectic}, if $w$ has type $\epsilon$ then $T_w$ is contained in a copy of $O^\epsilon_{2n}(q)$ but no copy of $O^{-\epsilon}_{2n}(q)$. So each separable element is contained in a subgroup isomorphic to $O^+_{2n}(q)$ or $O^-_{2n}(q)$ but not both.
	
	\subsection{Unitary groups}
	
	Let $X=SL_n(K)$ act on a $K$-vector space $E$ with basis $\{v_1, ..., v_n\}$. Take the maximal torus $R$ to be the group of diagonal matrices in $X$. The Weyl group $W\cong S_n$ is the set of permutation matrices. Take the Steinberg endomorphism $\sigma: (a_{ij}) \mapsto \left((a_{ij}^q)^{-\tr}\right)^\tau$ where $\tau$ is the matrix with $1$s on the anti-diagonal and $0$s elsewhere. Then $G=SU_n(q)$ which we view as a subgroup of $SL_n(q^2)$ which acts on the $\F_{q^2}$-space $V$ with basis $\{v_1, ..., v_n\}$ and preserves the non-singular conjugate-symmetric sesquilinear form
	\[\left(\, , \,\right):V \times V \rightarrow \F_{q^2}, \quad \left(\sum_{i=1}^n x_iv_i, \sum_{i=1}^n y_iv_i\right)=\sum_{i=1}^n x_iy_{n-i}^q.\] For each $d$, there is exactly one $G$-orbit of non-singular spaces of dimension $d$, and exactly one $G$-orbit of totally singular spaces of dimension $d$.
	
	The action of $\sigma$ on $W$ is non-trivial, and $W.\langle \sigma \rangle\cong W\times C_2$. Since this is a direct product, the $G$-conjugacy classes of $\sigma$-stable maximal tori correspond to partitions of $n$. We will work with the natural action of $W.\langle \sigma \rangle$ on the set $M := \{1, ..., n\}$.
	
	\begin{proposition}\mylabel{unitary}{Proposition~\thetheorem{}}
		If $w\in W$ is associated with the partition $n=n_1 + \cdots + n_m$ then $V$ decomposes under the action of $T_w$ as $V=V_1 \oplus \cdots \oplus V_m$ where $\dim V_i=n_i$ and $T_w$ stabilises the non-singular space $V_i$ for all $i$. If  $n_i$ is odd then $T_w$ acts irreducibly on $V_i$, whereas if $n_i$ is even then $T_w$ acts irreducibly on two spaces of dimension $n_i/2$ which sum to $V_i$.
	\end{proposition}
	\begin{proof}
		By \ref{correspondence}, specifically (\ref{compute Tw}), we know \[T_w\cong D_1 \times \cdots \times D_m\] where $D_i$ stabilises the space $V_i=V_i^+ \oplus V_i^-$ where $V_i^+$ has basis $\{v_j \mid j \mbox{ is even}, 1+\sum_{k=1}^{i-1}n_k \le j \le \sum_{k=1}^{i}n_k\}$ and $V_i^-$ has basis $\{v_{j} \mid j \mbox{ is odd}, 1+\sum_{k=1}^{i-1}n_k \le j \le \sum_{k=1}^{i}n_k\}$. The spaces $V_i^+, V_i^-$ are totally singular while $V_i$ is non-singular.
		
		Over the algebraically closed field we may write $D_i$ as \[D_i\cong \left\{\diag\left(t, t^{(-q)}, ..., t^{(-q)^{n_i-1}}\right) : t^{(-q)^{n_i}-1}=1\right\}\] so if $n_i$ is odd then $D_i$ acts irreducibly on $V_i$, whereas if $n_i$ is even then $D_i$ acts irreducibly on both $V_i^+$ and $V_i^-$.
	\end{proof}
	
	\subsection{Odd-dimensional orthogonal groups}

	Let $V$ be the $\F_q$-vector space with basis $\mathcal{B}:=\{v_0, v_1, ..., v_n, v_{-n}, ..., v_{-1}\}$ and define the quadratic form \[Q: V \rightarrow \F_q, \quad Q\left(\sum_{i=-n}^nx_iv_i\right)=x_0^2+\sum_{i=1}^n x_ix_{-i}.\] In odd characteristic $Q$ is non-singular, whereas in even characteristic $\rad(Q)$ is $1$-dimensional. The isometry group of $Q$ is $O_{2n+1}(q)$ which is isomorphic to $Sp_{2n}(q)$ in even characteristic. Since this case is treated in Section~\ref{symplectic section}, we will assume that $q$ is odd when working with odd-dimensional orthogonal groups.
	
	Let $X=Spin_{2n+1}(K)$ with $q$ odd and consider the action on the $K$-vector space $E$ with the above basis $\mathcal{B}$ as the image of the surjection $\phi: X \rightarrow SO_{2n+1}(K)$ defined in (\ref{surjection}) which preserves the extension of the form $Q$ to $E$. Take the Steinberg endomorphism $\sigma: (a_{ij}) \mapsto (a_{ij}^q)$. The fixed points form a group isomorphic to $Spin_{2n+1}(q)$ which is a double cover of $\Omega_{2n+1}(q)$, and we will work with the quotient group $G=\Omega_{2n+1}(q)$. The group $G$ acts on the $\F_q$-vector space $V$ described in the previous paragraph, preserving the quadratic form $Q$.
	
	We defined the type of a proper non-singular subspace in \ref{spacesdef}. Note that if $U_1$ and $U_2$ are non-singular even-dimensional complementary subspaces of $V$ then $U=U_1 \perp U_2$ has plus type if $U_1$ and $U_2$ have the same type, and minus type otherwise. For each $d$ there are exactly two orbits of proper $d$-dimensional non-singular subspaces of $V$ under the action of $G$: the spaces of plus type and those of minus type.
	
	Using the group of diagonal matrices in $X$ as our maximal torus $R$ we can write the Weyl group as a set of permutation matrices in $S_{2n}$ isomorphic to $B_n$, and so $\sigma$ acts trivially on $W$. Therefore $W.\langle\sigma\rangle \cong W\cong B_n$ which we will view acting on $M:=\{\pm 1, ..., \pm n\}$ in the natural way.
	
	\begin{proposition}\mylabel{odd orthogonal}{Proposition~\thetheorem{}}
		If $w \in W$ is in the conjugacy class associated to $(n_1^{\epsilon_1}, \dots, n_m^{\epsilon_m})$ then $V$ decomposes under the action of $T_w$ as $V=V_0 \oplus V_1 \oplus \cdots \oplus V_m$ where $\dim V_i = 2n_i$, $\dim V_0 = 1$ and $V_0$ is of $\prod_d \epsilon_d$ type. If $\epsilon_i=1$ then $V_i$ is non-singular of plus type and is a sum of two totally singular $n_i$-dimensional spaces upon which $T_w$ acts irreducibly. If $\epsilon_i=-1$ then $V_i$ is non-singular of minus type and $T_w$ acts on it irreducibly.
	\end{proposition}
	\begin{proof}
		By \ref{correspondence} --- see in particular (\ref{compute Tw}) --- we know $T_w\cong D_0 \times D_1^{\epsilon_1} \times \cdots \times D_m^{\epsilon_m}$ where $D_i^{\epsilon_i}$ stabilises the space $V_i=V_i^+ \oplus V_i^-$ where $V_i^+$ has basis $\{v_{j} \mid 1+\sum_{k=1}^{i-1}n_k \le j \le \sum_{k=1}^{i}n_k\}$ and $V_i^-$ has basis $\{v_{-j} \mid 1+\sum_{k=1}^{i-1}n_k \le j \le \sum_{k=1}^{i}n_k\}$. The spaces $V_i^+, V_i^-$ are totally singular while $V_i$ is non-singular.
		
		Over the field $K$, the group $D_i^+$ may be written\[D_i^+\cong \left\{\diag(t, t^q, ..., t^{q^{n_i-1}}, t^{-q^{n_i-1}},..., t^{-q}, t^{-1}) : t^{q^{n_i}-1}=1\right\}\]which acts irreducibly on the $n_i$-dimensional totally singular spaces $V_i^+$ and $V_i^-$, thereby stabilising $V_i$ which has plus type under the action of $T_w$.
		
		Over the field $K$, \[D_i^-\cong \left\{\diag(t^{-1}, t^q, ..., t^{q^{n_i-1}}, t^{-q^{n_i-1}},..., t^{-q}, t) : t^{q^{n_i}+1}=1\right\}\] which acts irreducibly on the $2n_i$-dimensional non-singular space $V_i$. Since there is no $n_i$-dimensional totally singular subspace under the action of $T_w$ in this case, $V_i$ is of minus type.
	
		By definition, the space $V_0$ which is centralised by $T_w$ has the same type as $w$.
	\end{proof}

	\subsection{Even-dimensional orthogonal groups}
	Fix a basis $\{v_1, ..., v_n, v_{-n}, ..., v_{-1}\}$ for the $K$-vector space $E=K^{2n}$ and let $X=Spin_{2n}(K)$. In odd characteristic let $X$ act on $E$ via the surjection $\phi: X \rightarrow SO_{2n}(K)$ defined in (\ref{surjection}). The action of $X$ on $E$ preserves the quadratic form \[Q:E \rightarrow K, \quad Q\left(\sum_{i=1}^nx_iv_i+x_{-i}v_{-i}\right)= \sum_{i=1}^nx_ix_{-i}.\] Using the group of diagonal matrices in $X$ as our maximal torus $R$ we can write the Weyl group as a set of permutation matrices in $S_{2n}$ isomorphic to $B_n^+$.
	
	The set of fixed points of the Steinberg endomorphism $\sigma_+: (a_{ij}) \mapsto (a_{ij}^q)$ in even characteristic is $G=\Omega_{2n}^+(q)$, and in odd characteristic forms a group isomorphic to $Spin_{2n}^+(q)$, a double cover of $\Omega_{2n}^+(q)$, and we take the quotient $G=\Omega^+_{2n}(q)$. Let $V$ be the $\F_q$-vector space with basis $\{v_1, ..., v_n, v_{-n}, ..., v_{-1}\}$. The elements of $G$ act on $V$ preserving the non-singular quadratic form of plus type \[Q_+:V \rightarrow \F_q, \quad Q_+\left(\sum_{i=1}^nx_iv_i+x_{-i}v_{-i}\right)= \sum_{i=1}^nx_ix_{-i}.\] Since $\sigma_+$ acts trivially on $W$ we have $W.\langle \sigma_+ \rangle\cong W\cong B_n^+$ which we will view acting on $M:=\{\pm 1, ..., \pm n\}$ in the natural way. The conjugacy classes correspond to signed cycle types $(n_1^{\epsilon_1}, \dots, n_r^{\epsilon_r})$ where $\prod_i\epsilon_i=1$.
	
	Consider the Steinberg endomorphism $\sigma_-: (a_{ij}) \mapsto (a_{ij}^q)^\tau$ where $\tau$ is a matrix representing the transposition $(n, -n)$ in $B_n$.  The set of fixed points of $\sigma_-$ in characteristic $2$ is $G=\Omega_{2n}^-(q)$, and in odd characteristic $p$ we have a group isomorphic to $Spin_{2n}^-(q)$ which is a double cover of $\Omega^-_{2n}(q)$, and we take the quotient $G=\Omega^-_{2n}(q)$. The elements of $G$ act on $V$ preserving the non-singular quadratic form of minus type \[Q_-:V \rightarrow \F_q, \quad Q_-\left(\sum_{i=1}^nx_iv_i+x_{-i}v_{-i}\right)= x_n^2+\mu x_{-n}^2+\sum_{i=1}^nx_ix_{-i}
	\] where $x^2+x+\mu$ is irreducible over $\F_q$. The endomorphism induced by $\sigma_-$ on $W$ has order $2$ and $W.\langle \sigma_- \rangle\cong B_n$ which acts on $M:=\{\pm 1, ..., \pm n\}$ in the natural way. The coset $W\sigma_-$ is the set of negative elements of $B_n$, and the $W$-conjugacy classes correspond to signed cycle types $(n_1^{\epsilon_1}, \dots, n_r^{\epsilon_r})$ where $\prod_i\epsilon_i=-1$.
	
	For $G$ of plus or minus type, for each $d$ there is exactly one $G$-orbit of proper $d$-dimensional non-singular subspaces of plus type of $V$, and exactly one $G$-orbit of proper $d$-dimensional non-singular subspaces of minus type of $V$.
	
	\begin{proposition}\mylabel{even orthogonal}{Proposition~\thetheorem{}}
		For $\sigma \in \{\sigma_+, \sigma_-\}$ if $w.\sigma \in W.\langle \sigma \rangle$ is in the $W$-conjugacy class associated to $(n_1^{\epsilon_1}, \dots, n_m^{\epsilon_m})$ then $V$ decomposes as $V= V_1 \oplus \cdots \oplus V_m$ where $\dim V_i = 2n_i$. If $\epsilon_i=1$ then $V_i$ is non-singular of plus type and is a sum of two totally singular $n_i$-dimensional spaces upon which $T_w$ acts irreducibly. If $\epsilon_i=-1$ then $V_i$ is non-singular of minus type and $T_w$ acts on it irreducibly.
	\end{proposition}
	\begin{proof}
		If $G=\Omega_{2n}^-(q)$ then $w$ has minus type and so we may relabel the $n_i$ so that $\epsilon_m=-1$. By \ref{correspondence} --- see in particular (\ref{compute Tw}) --- we know $T_w\cong D_1^{\epsilon_1} \times \cdots \times D_m^{\epsilon_m}$ where $D_i^{\epsilon_i}$ stabilises the space $V_i=V_i^+ \oplus V_i^-$ where $V_i^+$ has basis $\{v_{j} \mid 1+\sum_{k=1}^{i-1}n_k \le j \le \sum_{k=1}^{i}n_k\}$ and $V_i^-$ has basis $\{v_{-j} \mid 1+\sum_{k=1}^{i-1}n_k \le j \le \sum_{k=1}^{i}n_k\}$. The space $V_i$ is non-singular. The spaces $V_i^+, V_i^-$ are totally singular unless $G=\Omega_{2n}^-(q)$ and $i=m$ in which case $v_n\in V_m^+$ and $v_{-n} \in V_m^-$, and by assumption $\epsilon_i=-1$. Over the algebraic closure $K$, we may write \[D_i^+\cong \left\{\diag(t, t^q, ..., t^{q^{n_i-1}}, t^{-q^{n_i-1}},..., t^{-q}, t^{-1}) : t^{q^{n_i}-1}=1\right\}\] which acts irreducibly on $V_i^+, V_i^-$, so $V_i$ has plus type and \[D_i^{-}\cong \left\{\diag(t^{-1}, t^q, ..., t^{q^{n_i-1}}, t^{-q^{n_i-1}},..., t^{-q}, t) : t^{q^{n_i}+1}=1\right\}\] which acts irreducibly on $V_i$ and does not stabilise an $n_i$-dimensional totally singular subspace, so $V_i$ has minus type.
	\end{proof}

\begin{corollary}\mylabel{sets and spaces}{Corollary~\thetheorem{}}
	Take $w \in W$ and $g$ a separable element in some conjugate of $T_w$.
	\begin{enumerate}[(a)]
		\item If $G=SL_n(q)$ then $\Inv(g)=\Inv(w)$.
		\item If $G=Sp_{2n}(q)$ then $\Inv_\mathcal{N}(g)=\Inv_\mathcal{N}(w)$.
		\item If $G=SU_{n}(q)$ then $\Inv_\mathcal{N}(g)=\Inv(w)$.
		\item If $G=\Omega_{2n+1}(q)$ then for positive $g$ we get $\Inv_{\mathcal{N}^+}(g)=\Inv_{\mathcal{N}^+}(w) \cup \{1\}$ and $\Inv_{\mathcal{N}^-}(g)=\Inv_{\mathcal{N}^-}(w)$. If $g$ is negative then $\Inv_{\mathcal{N}^+}(g)=\Inv_{\mathcal{N}^+}(w)$ and $\Inv_{\mathcal{N}^-}(g)=\Inv_{\mathcal{N}^-}(w)  \cup \{1\}$.
		\item If $G=\Omega^\pm_{2n}(q)$ then $\Inv_{\mathcal{N}^+}(g)=\Inv_{\mathcal{N}^+}(w)$ and $\Inv_{\mathcal{N}^-}(g)=\Inv_{\mathcal{N}^-}(w)$.
	\end{enumerate}
\end{corollary}
\begin{proof}
	First assume $G=SL_n(q)$. If $w$ corresponds to the partition $n_1+\cdots + n_m=n$ then it follows that the invariant sets are all the possible unions of the invariant sets of each cycle. So $\Inv(w)$ is the set of subsums \[\Inv(w)=\left\{\sum_{j \in J}n_j \,\Big|\, \emptyset \ne J \subsetneq \{1, ..., m\}\right\}.\] By \ref{torus element}, $g$ acts in the same way as $T_w$ and so by \ref{linear}, and since the invariant spaces are direct sums of irreducible invariant spaces, $\Inv(g)=\Inv(w)$.
	
	Now assume $G=Sp_{2n}(q)$ and take $w$ with signed cycle type $(n_1^{\epsilon_1}, \dots, n_m^{\epsilon_m})$. Then, arguing in the same way as above, \[\Inv_\mathcal{N}(w) = \left\{\sum_{j \in J}2n_j \,\Big|\, \emptyset \ne J \subsetneq \{1, ..., m\}\right\}\] and by \ref{symplectic} we can compute	$\Inv_\mathcal{N}(g)=\Inv_\mathcal{N}(w)$.
	
	For $G=SU_n(q)$ take $w$ corresponding to the partition $n_1+\cdots + n_m=n$ so that \[\Inv(w)=\left\{\sum_{j \in J}n_j \,\Big|\, \emptyset \ne J \subsetneq \{1, ..., m\}\right\}\] and by \ref{unitary} we get $\Inv_\mathcal{N}(g)=\Inv(w)$.
	
	Take $G=\Omega_{2n+1}(q)$ with $w$ of signed cycle type $(n_1^{\epsilon_1}, \dots, n_m^{\epsilon_m})$. Arguing as above, \[\Inv_\mathcal{N^+}(w)=\left\{\sum_{j \in J}2n_j \,\Big|\, \emptyset \ne J \subsetneq \{1\le i \le m : \epsilon_i=+\}\right\}\] and by \ref{odd orthogonal} \[\Inv_\mathcal{N^+}(g)=\begin{cases}
		\Inv_\mathcal{N^+}(w) & w \mbox{ of minus type}\\
		\Inv_\mathcal{N^+}(w) \cup \{1\} & w \mbox{ of plus type.}
	\end{cases}\]We may use the same argument for $\mathcal{N}^-$.

	Finally let $G=\Omega_{2n}^\epsilon(q)$ and take $w$ with signed cycle type $(n_1^{\epsilon_1}, \dots, n_m^{\epsilon_m})$. Arguing as above and by \ref{even orthogonal} \[\Inv_\mathcal{N^+}(g)=\left\{\sum_{j \in J}2n_j \,\Big|\, \emptyset \ne J \subsetneq \{1\le i \le m : \epsilon_i=+\}\right\}= \Inv_\mathcal{N^+}(w)\] and similarly for $\mathcal{N}^-$.
\end{proof}

	\section{Quasi-invariant spaces and sets}\label{quasi}
	
	Recall that $s(G)$ is the set of separable elements in $G$. If $g \in s(G)$ then we can take $w \in W$ such that $g$ is contained in some conjugate of $T_{w}$. In fact, the following discussion shows that the proportion of elements in $G$ which are separable and in some conjugate of $T_w$ can be approximated by the proportion of elements in $W$ which are conjugate to $w$.
	
	Pick a set of representatives $C$ of the $W$-conjugacy classes in the coset $W\sigma$ of $W.\langle\sigma\rangle$, and for $w \in C$ define $w^G$ to be the set of separable elements of $G$ contained in some conjugate of $T_w$, mirroring the notation $(w\sigma)^W$ for the $W$-conjugacy class in $W.\langle\sigma\rangle$. For $D \subseteq C^k$ a set of $k$-tuples of $W$-conjugacy classes in $W\sigma$ define the following disjoint unions \begin{equation}\label{disjoint unions}D^W=\bigcup_{(d_1, ..., d_k) \in D}\prod_{i=1}^k (d_i\sigma)^W, \quad D^G=\bigcup_{(d_1, ..., d_k) \in D}\prod_{i=1}^k d_i^G\end{equation} so that $D^W$ is the set of all conjugates of tuples in $D$, and $D^G$ is the set of all conjugates of tuples in $s(G)$ which correspond to tuples in $D$.
	
	Recall the definition of $\epsilon(q)$ (\ref{epsilon}) which is an upper bound for the proportion of elements in $G$ which are not separable. The following proposition will be very useful throughout the rest of the paper.
	
	\begin{proposition}\mylabel{proportion}{Proposition~\thetheorem{}}
		For any $D \subseteq C^k$,
		\[
		\frac{|D^W|}{|W|^k} -1+\left(1- \epsilon(q)\right)^k \le \frac{|D^G|}{|G|^k}\le \frac{|D^W|}{|W|^k}
		\]
	\end{proposition}
	\begin{proof}
		Recall from \ref{most sep} that \[
		\frac{\left |s(G)\right |}{\left| G \right |}\ge 1-\epsilon(q).\] Fulman and Guralnick \cite[Section~5]{FulmanGuralnick} show that \[\frac{|w^G|}{|G|}\le \frac{|w\sigma^W|}{|W|}\] for all $w \in C$ and by the fact that the unions in (\ref{disjoint unions}) are disjoint, it follows that, for any $D \subseteq C^k$, \[\frac{|D^G|}{|G|^k}\le \frac{|D^W|}{|W|^k}.\] In particular we also get \[\frac{\left\vert\left(C^k\setminus D\right)^G\right\vert}{|G|^k}\le \frac{\left\vert\left(C^k\setminus D\right)^W\right\vert}{|W|^k}\] and since $\left\vert\left(C^k\setminus D\right)^G\right\vert+\left\vert D^G\right\vert=|s(G)|^k$ and $\left\vert\left(C^k\setminus D\right)^W\right\vert+\left\vert D^W\right\vert=|W|^k$ the result follows.
	\end{proof}

	Now we will describe the proportion of subsets of size $k$ in $G$ that are made up of separable elements which invariably generate an irreducible subgroup of $G$ in terms of the spaces they stabilise.
	
	Recall that $G$ is a finite classical group from the list (\ref{classical groups}) acting on $V$ with Weyl group $W$ acting on $M$ as described in Section~\ref{invariant}. That is, for $G=SL_n(q), SU_n(q)$ we take $M=\{1, ..., n\}$ and for $G=Sp_{2n}(q), \Omega_n^\epsilon(q)$ we take $M=\{\pm 1, ..., \pm n\}$.
	
	\begin{definition}\label{quasi invariant}
		For a subset $S=\{g_1, ..., g_k\} \subseteq G$, a subspace $U \le V$ will be called $S$-quasi-invariant if there exist $h_1, ..., h_k \in G$ such that $g_i^{h_i}(U)=U$ for $i=1, ..., k$. For convenience we will write $g$-quasi-invariant when we mean $\{g\}$-quasi-invariant.
	\end{definition}

	The following proposition is the main result of this section.
	
	\begin{proposition}\mylabel{invariant spaces}{Proposition~\thetheorem{}}
		Pick $g_1, ..., g_k$ independently uniformly at random from $G$ and let $S=\{g_1, ..., g_k\}$. There are functions $f_G(n, q)$ with values in $[0, 1]$ such that, for fixed $q$, $\lim_{n \rightarrow \infty}f_G(n,q)=0$ and
		\begin{enumerate}[(a)]
			\item $\PI(G, k, s(G), \I)=\Prob(\{\cap \Inv(g_i)=\emptyset\}\cap \{S\subseteq s(G)\})$ for $G=SL_n(q)$,
			\item $\PI(G, k, s(G), \I)=\Prob(\{\cap \Inv_\mathcal{N}(g_i)=\emptyset\}\cap \{S\subseteq s(G)\})+f_G(n,q)$ for $G=Sp_{2n}(q)$ or $SU_{n}(q)$,
			\item $\PI(G, k, s(G), \I)=\Prob(\{\cap \Inv_\mathcal{N^+}(g_i)=\emptyset\} \cap \{\cap \Inv_\mathcal{N^-}(g_i)=\emptyset\}\cap \{S\subseteq s(G)\})+f_G(n,q)$ for $G=\Omega_{2n+1}(q), \Omega^+_{2n}(q)$ or $\Omega^-_{2n}(q)$.
		\end{enumerate}
	\end{proposition}

	Denote by $\mathcal{U}(S)$ the set of all non-zero proper $S$-quasi-invariant subspaces of $V$ and say $\mathcal{U}_\mathcal{N}(S)\subseteq \mathcal{U}(S)$ is the set of non-singular elements of $\mathcal{U}(S)$ and $\mathcal{U}_\mathcal{R}(S)\subseteq \mathcal{U}(S)$ is the set of totally singular elements of $\mathcal{U}(S)$. Say $S$ is bad if $S \subseteq s(G)$ and $\mathcal{U}_\mathcal{N}(S)=\emptyset$ and $\mathcal{U}(S)\ne\emptyset$. In order to prove \ref{invariant spaces} we need the following lemma which tells us that if $G\ne SL_n(q)$ then the proportion of bad sets vanishes as $n \rightarrow \infty$.
	\begin{lemma}\mylabel{nondegenerate}{Lemma~\thetheorem{}}
		If $G\ne SL_n(q)$ and \[f_G(n,q):= \frac{\{S \subseteq s(G) : |S|=k \mbox{ and }S\mbox{ is bad}\}}{|G|^k}\] then, fixing $q$, $\lim_{n \rightarrow \infty}f_G(n,q)=0.$
	\end{lemma} 
	\begin{proof}
		Assuming $\mathcal{U}(S)$ is non-empty, take $U \in \mathcal{U}(S)$ of minimal dimension. Then either $\rad(U)=0$ or $\rad(U) \in \mathcal{U}(S)$, in which case $\rad(U)=U$ by minimality. So if $\mathcal{U}(S)\ne\emptyset$ then either $\mathcal{U}_\mathcal{N}(S)\ne\emptyset$ or $\mathcal{U}_\mathcal{R}(S)\ne\emptyset$. Now $S$ is bad if and only if $\mathcal{U}_\mathcal{N}(S)=\emptyset$ and $\mathcal{U}_\mathcal{R}(S)\ne\emptyset$.
		
		Section~\ref{invariant} describes the situations giving rise to bad sets $S$. There are no bad sets when $G=\Omega_{2n+1}(q)$ in odd characteristic, when $G=\Omega^-_{2n}(q)$, or when $G=SU_{2n+1}(q)$. For example, if $S$ is a bad set for $G=\Omega^-_{2n}(q)$ then $\mathcal{U}(S)$ consists only of two totally singular spaces of dimension $n$, which can only occur when each cycle of each $w_i$ has plus type.
		
		For $G=Sp_{2n}(q)$ or $SU_{2n}(q)$ (see \ref{symplectic} and \ref{unitary}), if $S$ is bad then all cycles in all $w_i$ have even length. By \ref{proportion} and \cite[Lemma 4.2]{FulmanGuralnick}, the contribution of this case tends to $0$ as $n$ tends to infinity.
		
		For $G=\Omega_{2n}^+(q)$, \ref{even orthogonal} shows that if $S$ is bad then all cycles in all $w_i$ are positive. By \ref{proportion} and \cite[Theorem 4.4]{FulmanGuralnick}, the contribution of this case tends to $0$ as $n$ tends to infinity.
	\end{proof}

	We have defined $f_G$ when $G\ne SL_n(q)$. For $G=SL_n(q)$ we may set $f_G(n, q):=0$ for all $n, q$.
		
	\begin{proof}[Proof of \ref{invariant spaces}]
		The set $S$ fails to invariably generate an irreducible subgroup of $G$ if and only if there is a non-zero proper $S$-quasi-invariant subspace $U<V$. First assume $G=SL_n(q)$. Since each $d$-dimensional space is in the same $G$-orbit, if $g\in G$ stabilises a space of dimension $d$ then every $d$-dimensional space is $g$-quasi-invariant. So $S$ fails to invariably generate an irreducible subgroup of $G$ if and only if $d \in \cap_i \Inv{g_i}$ for some $0<d<n$.

		For the remaining cases, there is a form on $V$. By \ref{nondegenerate}, the contribution of bad sets $S$ is $f_G(n,q)$ which vanishes as $n$ tends to infinity. So the probability that $S\subseteq s(G)$ and $S$ invariably generates an irreducible subgroup of $G$ is $\beta + f_G(n,q)$ where $\beta$ is the probability that $S\subseteq s(G)$ and there is no non-zero proper non-singular $S$-quasi-invariant subspace.
		
		Now assume $G=Sp_{2n}(q)$ or $SU_n(q)$. Arguing as above, since each non-singular $d$-dimensional space is in the same $G$-orbit, $\beta=\Prob(\{\cap \Inv_\mathcal{N}(g_i)=\emptyset\}\cap \{S\subseteq s(G)\})$.
		
		Finally take $G=\Omega_{2n+1}(q), \Omega^+_{2n}(q)$ or $\Omega^-_{2n}(q)$. There are exactly two orbits of non-singular spaces of dimension $d$ under the action of $G$: the spaces of plus type and those of minus type, so arguing as before $\beta=\Prob(\{\cap \Inv_\mathcal{N^+}(g_i)=\emptyset\} \cap \{\cap \Inv_\mathcal{N^-}(g_i)=\emptyset\}\cap \{S\subseteq s(G)\})$.
	\end{proof}
	
	\section{Distribution of separable elements and Weyl group elements}\label{irred sep}

	Pick $g_1, ..., g_k$ independently uniformly at random from $G$ and let $S:=\{g_1, ..., g_k\}$. If $S \subseteq s(G)$ then we may pick $w_1, \ldots, w_k \in W$ such that $g_i$ is contained in some conjugate of $T_{w_i}$.	Now select $z_1, ..., z_k \in W$ uniformly at random. Fulman and Guralnick \cite[Section~5]{FulmanGuralnick} show that the $w_i$ and $z_i$ have similar distributions, as discussed in the proof of \ref{proportion}.
	
	For any set $\Property$ of subsets of $M$, setting \[D=\left\{(d_1, ..., d_k) \mid \cap \Inv_\Property(d_i)=\emptyset\right\}\subseteq C^k\] we have \[\Prob\left(\cap\Inv_\Property(z_i)=\emptyset\right)= \frac{|D^W|}{|W|^k}, \quad \Prob\left(\{S \subseteq s(G)\}\cap \{\cap\Inv_\Property(w_i)=\emptyset\}\right)=\frac{|D^G|}{|G|^k}.\]
	
	Define $\delta(q,k)=1-\left(1- \epsilon(q)\right)^k>0$ for all $k$ and $q$ and note that $\delta(q,k)=O(1/q)$ as $q \rightarrow \infty$ for any $k$ and $n$. We will need the following corollary of \ref{proportion}.
	
	\begin{corollary}\mylabel{proportion cor}{Corollary~\thetheorem{}}
		For any set $\Property$ of subsets of $M$, \[\Prob(\cap\Inv_\Property(z_i)=\emptyset) - \delta(q,k) \le \Prob(\{S\subseteq s(G)\}\cap\{\cap\Inv_\Property(w_i)=\emptyset\}) \le \Prob(\cap\Inv_\Property(z_i)=\emptyset).\]
	\end{corollary}
	
	To conclude this section, we describe the relationship between the invariable generation properties of a finite classical group and its Weyl group when we restrict to separable elements.
	
	Recall that $\PI(G,k,H,\mathcal{A})$ is the probability that the set $S=\{g_1, ..., g_k\}$ of independent random elements of $G$ is contained in $H$ and invariably generates a member of $\mathcal{A}$. We also defined $\mathcal{T}$ to be the set of proper transitive subgroups of a permutation group, and $\mathcal{I}$ to be the set of irreducible subgroups of $G$.

	\begin{theorem}\mylabel{probably same}{Theorem~\thetheorem{}}
	For $k \ge 2$ \[\left(1-2^{1-k}\right)\PI(S_n, k, \T)-3\delta(q,k)-\alpha(n, k)\le \PI(G, k).\] In addition either \[\PI(G, k)\le\PI(S_n, k, \T)+\epsilon(q)^k+f_G(n, q)\] or $G$ is an orthogonal group and
	\[\PI(G, k)\le\Prob(\{\cap \Inv_\mathcal{N^+}(z_i)=\emptyset\} \cap \{\cap \Inv_\mathcal{N^-}(z_i)=\emptyset\})+\epsilon(q)^k+f_G(n, q).\]
\end{theorem}

	\begin{proof}
		First consider the case $G=SL_n(q)$. Applying \ref{invariant spaces} followed by \ref{sets and spaces},
		\[\PI(G, k, s(G), \I)=\Prob(\{S \subseteq s(G)\}\cap\{\cap \Inv(w_i)=\emptyset\})\]
		and so by \ref{proportion cor}
		\[\Prob(\cap \Inv(z_i)=\emptyset)-\delta(q,k)\le \PI(G, k, s(G), \I)\le\Prob(\cap \Inv(z_i)=\emptyset).\] Finally by \ref{rss irred} and the facts that $\I=\I'$ and $\PI(S_n, k, \T)=\Prob(\cap \Inv(z_i)=\emptyset)$ we see
		\[\PI(S_n, k, \T)-\delta(q,k)-\alpha(n, k)\le \PI(G, k)\le\PI(S_n, k, \T)+\epsilon(q)^k.\]
		
		The argument for $G=SU_n(q)$ is exactly the same.
		
		For the remaining cases let $\pi:W \rightarrow S_n$ be the canonical projection and note that \begin{align*}\Prob(\{\cap \Inv_\mathcal{N^+}(z_i)=\emptyset\} \cap \{\cap \Inv_\mathcal{N^-}(z_i)=\emptyset\})&\ge \Prob(\cap \Inv_{\mathcal{N}}(z_i)=\emptyset)\\&=\Prob(\{\cap \Inv(\pi(z_i))=\emptyset\})\\&= \PI(S_n, k, \T).\end{align*}
		
		For $G=Sp_{2n}(q)$ in odd characteristic, we may use the same argument as before. In even characteristic, we note $\I \ne \I'$, and we deal with this below.
		
		For $G=\Omega_{2n}^\pm(q)$ the same argument as above gives us \[\Prob(\{\cap \Inv_\mathcal{N^+}(z_i)=\emptyset\} \cap \{\cap \Inv_\mathcal{N^-}(z_i)=\emptyset\})-\delta(q,k)-\alpha(n, k)\le \PI(G, k)\]
		and also \[\PI(G, k)\le \Prob(\{\cap \Inv_\mathcal{N^+}(z_i)=\emptyset\} \cap \{\cap \Inv_\mathcal{N^-}(z_i)=\emptyset\})+\epsilon(q)^k+f_G(n,q).\]
		
		In the remaining cases, let $A$ be the event that the random elements $g_i$ in $G$ are all separable and the resulting $w_i$ do not all have the same sign. By \ref{proportion cor}, $\Prob(A)\ge 1-2^{1-k}-\delta(q,k)$ since the probability that the $z_i$ do not all have the same sign is $1-2^{1-k}$.
		
		For $G=Sp_{2n}(q)$ in even characteristic, $\I \ne \I'$. For $i=1, ..., k$, if $w_i$ is of $\epsilon$ type where $\epsilon\in \{+, -\}$ then $g_i$ is contained in some subgroup isomorphic to $O^\epsilon_{2n}(q)$, and no subgroup isomorphic to $O^{-\epsilon}_{2n}(q)$. Therefore by \ref{invariant spaces} \[\PI(G,k,s(G),\I')=\Prob(A \cap \{S\subseteq s(G)\}\cap\{\cap \Inv_{\mathcal{N}}(w_i)=\emptyset\})+f_{G}(n, q).\]
		
		Applying \ref{proportion cor}, noting that $\Inv_{\mathcal{N}}(w_i)$ is independent of the sign of $w_i$ by \ref{symplectic}, we get \begin{align*}\PI(G,k,s(G),\I')&\ge\Prob(A)(\Prob(\cap \Inv_{\mathcal{N}}(z_i)=\emptyset)-\delta(k,q)+f_{G}(n, q))\\ &\ge (1-2^{1-k})\Prob(\cap \Inv_{\mathcal{N}}(z_i)=\emptyset)-3\delta(k,q).\end{align*}
		
		For $G=\Omega_{2n+1}(q)$, the argument for even dimensional orthogonal groups gives us the upper bound on $\PI(G, k)$. For the lower bound, note that elements which invariably generate cannot all share the same sign. Then arguing as above \[\PI(G, k)\ge (1-2^{1-k})\Prob(\{\cap \Inv_\mathcal{N^+}(z_i)=\emptyset\} \cap \{\cap \Inv_\mathcal{N^-}(z_i)=\emptyset\})-3\delta(q,k)-\alpha(n, k)\]
		which finishes the proof.
	\end{proof}
	
	\section{Invariable generation of Weyl groups}\label{weyl}
	
	\begin{lemma}\mylabel{generalEFG}{Lemma~\thetheorem{}}
		For $k=3$ and $W=B_n$ or $B_n^+$, \[\lim_{n \rightarrow \infty}\Prob(\{\cap \Inv_\mathcal{N^+}(z_i)=\emptyset\} \cap \{\cap \Inv_\mathcal{N^-}(z_i)=\emptyset\})=0.\]
	\end{lemma}
	
	The author thanks Sean Eberhard for suggesting the following proof which is a modification of the proof of \ref{efg}(b).
	
	\begin{proof}
		Recall that the $z_i$ are random elements of $W$ and let $z_i'$ for $i=1, 2, 3$ be elements selected uniformly at random from $S_n$.
		
		We define versions of $\Inv$ which only consider cycles whose length is in a set $I\subseteq \mathbb{R}$: if $w'\in S_n$ corresponds to the partition $n_1 + \cdots + n_m=n$ and if $w\in W$ has signed cycle type $(n_1^{\epsilon_1}, \dots, n_m^{\epsilon_m})$ then define \[\Inv^I(w)=\Inv^I(w')=\left\{\sum_{j \in J}n_j \,\Big|\, \emptyset \ne J \subseteq \{i : n_i \in I\}\right\}\]
		along with the following variations
		\begin{align*}\Inv^I_\mathcal{N^+}(w)&=\left\{\sum_{j \in J}2n_j \,\Big|\, \emptyset \ne J \subseteq \{i : n_i \in I, \epsilon_i=+\}\right\},\\ \Inv^I_\mathcal{N^-}(w)&=\left\{\sum_{j \in J}2n_j \,\Big|\, \emptyset \ne J \subseteq \{i : n_i \in I, \epsilon_i=-\}\right\}.\end{align*}
		
		Eberhard, Ford and Green \cite[Proposition~3.8]{EFG2} show that there are constants $d>0$ and \[\beta=1-2(3 \log 2)^{-1}-0.02\] such that for any interval $I=(j^\beta, 60j]$ we have $\Prob(\cap_{i=1}^3 \Inv^I(z_i')\ne \emptyset)\ge d$. Summing over infinitely many disjoint intervals $I$ they show that $\lim_{n \rightarrow \infty}\Prob(\cap_{i=1}^3 \Inv(z_i'))=1.$ For $W=B_n$ or $B_n^+$, we may pick the $z_i$ by picking elements $z_i' \in S_n$ and then independently assigning signs to each cycle uniformly at random.
		
		If $l \in \cap_{i=1}^3 \Inv^I(z_i)$, then $l$ is contained in $\cap_{i=1}^3 \Inv^I_{\mathcal{N}^+}(z_i)$ or $\cap_{i=1}^3 \Inv^I_{\mathcal{N}^-}(z_i)$ with probability $\frac{1}{4}$, and so for any interval $I=(j^\beta, 60j]$, we have \[\Prob(\{\cap \Inv^I_\mathcal{N^+}(z_i)\ne\emptyset\} \cup \{\cap \Inv^I_\mathcal{N^-}(z_i)\ne\emptyset\})\ge \frac{d}{4}.\] Now, following \cite{EFG2}, we can sum over the same infinite set of intervals to see that \[\lim_{n \rightarrow \infty}\Prob(\{\cap \Inv_\mathcal{N^+}(z_i)\ne\emptyset\} \cup \{\cap \Inv_\mathcal{N^-}(z_i)\ne\emptyset\})=1.\]
	\end{proof}
	
	\section{Proof of \ref{main}}\label{proof}
	
	In this final section, we present a proof of \ref{main}. First, let us prove the result for a finite classical group $G$ from the list (\ref{classical groups}). By \ref{probably same}, for $k=4$ we have \[\PI(G, 4)\ge\frac{7}{8}\PI(S_n, 4, \T)-3\delta(q,4)-\alpha(n,4).\] By \ref{efg}(a), there is a constant $b>0$ such that $\PI(S_n, 4, \T)\ge b$ for all $n$. Therefore \[\PI(G, 4)\ge \frac{7}{8}b-3\delta(q,4)-\alpha(n,4).\] Since $\lim_{q \rightarrow \infty}\delta(q,4)=\lim_{n \rightarrow \infty}\alpha(n, 4)=0$ we have bounded $\PI(G, 4)$ away from $0$ for large enough $q$ and $n$. This works for groups of large rank, say $r>R$. For groups with smaller rank, note that by \ref{GM} there is an absolute constant $\epsilon_r$ such that $\PI(G,4)\ge \epsilon_r$ whenever $G$ is a finite group of Lie type of rank $r \le R$ and $q$ is large enough.
			
	By \ref{probably same}, if $G$ is not an orthogonal group, then \[\PI(G, 3)\le\PI(S_n, 3, \T)+\epsilon(q)^3+f_G(n, q)\] which tends to $0$ as $n \rightarrow \infty$ and $q \rightarrow \infty$ by \ref{efg}(b).
			
	For $G$ an orthogonal group,
	\[\PI(G, 3)\le\Prob(\{\cap \Inv_\mathcal{N^+}(z_i)=\emptyset\} \cap \{\cap \Inv_\mathcal{N^-}(z_i)=\emptyset\})+\epsilon(q)^3+f_G(n, q)\] which by \ref{generalEFG} tends to $0$ as $q\rightarrow\infty$ after $n \rightarrow \infty$.
	
	We have proved \ref{main} for the finite classical groups. Now it remains to extend this to all finite groups of Lie type. First note that the exceptional groups have bounded rank and so for part $(a)$ they are covered by \ref{GM} and we do not need to consider them when proving $(b)$.
	
	If $H$ is a finite group of Lie types $A, B, C$ or $D$ then $H=G/N$ where $G=X_\sigma$ for $X$ a simply connected linear algebraic group and $\sigma$ a Steinberg endomorphism, and $N\le Z(G)$. 
	%By \ref{most sep}, the proportion of separable elements in any coset of $N$ in $G$ is at least $1-\epsilon(q)$.
	Then every invariable generating set of $G$ projects onto one of $H$. For an invariable generating set $\{h_1, ..., h_k\}$ of $H$, a lift $\{h_1n_1, ..., h_kn_k\}$ with $n_i \in N$ invariably generates $G$ if and only if the set $\{n_1, ..., n_k\}$ invariably generates $N$. Therefore $\PI(G, k)=\PI(H,k)\PI(N,k)$. If we can show that $\PI(N,k)$ is bounded away from zero by an absolute constant $\gamma$ then we will have \begin{equation}\label{inequalities}
	\gamma\PI(H, k)\le\PI(G, k)\le \PI(H, k) \le \gamma^{-1}\PI(G, k).
	\end{equation}
		
	When $G\ne SL_n(q), SU_n(q)$, there are only finitely many options for $N$ and so $\PI(N,k)$ is bounded away from zero by an absolute constant.
		
	When $G=SL_n(q)$ let $d:=(n,q-1)$, and when $G=SU_n(q)$ let $d:=(n, q+1)$. Then $N$ is a subgroup of $\mathbb{Z}_d$ which can be viewed as the group of $d$th roots of unity. So $N$ is a cyclic group of order $e | d$. Letting $\zeta$ be the Riemann zeta function, $\PI(N,k)=\prod_{l | e}(1-1/l^k)>1/\zeta(k)$ where the product is over the set of prime divisors $l$ of $e$ (see, for example, \cite{pomerance2002expected}). So $\PI(N,3), \PI(N,4)$ are bounded away from zero by an absolute constant.
	
	So for $k=3, 4$ we have shown the inequalities in (\ref{inequalities}). Consequently, if the result holds for some quotient of $G$ by a central subgroup then it holds for all quotients of $G$ by central subgroups. This completes the proof for all finite groups of Lie type.
	\qed
	
	\begin{remark}If we know the value of the constant $b$ in \ref{efg}(a) then we can calculate the value of $Q$ in \ref{main}. Some experiments in GAP \cite{GAP4} lead us to guess that the value of $b$ could be close to $1/3$. If we have $b=1/3$ then we know $\PI(G, 4)\ge\frac{7}{24}-3\delta(q,4)-\alpha(n,4)$ so assuming that $n$ is large enough and $b=1/3$ we can say that $Q\le279$.
	\end{remark}
	\bibliographystyle{elsarticle-harv}
	\bibliography{../../bib/bib}
\end{document}